\documentclass[12pt,a4paper]{article}
\usepackage{amscd,amsmath,amssymb,amsfonts}

\numberwithin{equation}{section}

\newtheorem{thm}{Theorem}[section]

\newtheorem{lem}[thm]{Lemma}
\newtheorem{cor}[thm]{Corollary}

\newtheorem{claim}[thm]{Claim}

\newtheorem{rem}[thm]{Remark}
\def\Th#1.{\vskip 6pt \medbreak\noindent{\bf Theorem #1.}}

\newenvironment{pf}{\smallskip\noindent\emph{Proof.}}{Q.E.D.\bigbreak}

\begin{document}

\def\Spec{{\mathrm{Spec}}}
\def\Pic{{\mathrm{Pic}}}
\def\Ext{{\mathrm{Ext}}}
\def\NS{{\mathrm{NS}}}
\def\CH{{\mathrm{CH}}}
\def\deg{{\mathrm{deg}}}
\def\dim{{\mathrm{dim}}}
\def\codim{{\mathrm{codim}}}
\def\Coker{{\mathrm{Coker}}}
\def\ker{{\mathrm{ker}}}
\def\ch{{\mathrm{ch}}}
\def\Image{{\mathrm{Image}}}
\def\Aut{{\mathrm{Aut}}}
\def\Hom{{\mathrm{Hom}}}
\def\Proj{{\mathrm{Proj}}}
\def\Sym{{\mathrm{Sym}}}
\def\Image{{\mathrm{Image}}}
\def\Gal{{\mathrm{Gal}}}
\def\GL{{\mathrm{GL}}}
\def\End{{\mathrm{End}}}
\def\P{{\mathbb P}}
\def\C{{\mathbb C}}
\def\R{{\mathbb R}}
\def\Q{{\mathbb Q}}
\def\Z{{\mathbb Z}}
\def\F{{\mathbb F}}
\def\l{\ell}
\def\lra{\longrightarrow}
\def\ra{\rightarrow}
\def\hra{\hookrightarrow}
\def\ot{\otimes}
\def\op{\oplus}
\def\vg{\varGamma}
\def\O{{\cal{O}}}
\def\ol#1{\overline{#1}}
\def\wt#1{\widetilde{#1}}
\def\us#1#2{\underset{#1}{#2}}
\def\os#1#2{\overset{#1}{#2}}
\def\lim#1{\us{#1}{\varinjlim}}
\def\plim#1{\us{#1}{\varprojlim}}

%


\def\Gr{{\mathrm{Gr}}}
\def\rank{{\mathrm{rank}}}
\def\dlog{{\mathrm{dlog}}}
\def\Res{{\mathrm{Res}}}
\def\tor{{\mathrm{tor}}}
\def\reg{{\mathrm{reg}}}

\def\cJ{{\cal{J}}}
\def\cH{{\cal{H}}}
\def\cX{{\cal{X}}}
\def\K{{\cal{K}}}
\def\E{{\cal{E}}}
\def\G{{\mathbb{G}}}
\def\dR{{\mathrm{dR}}}
\def\DR{{\mathrm{DR}}}
\def\Zar{{\mathrm{Zar}}}
\def\prim{{\mathrm{prim}}}
\def\NS{{\mathrm{NS}}}
\def\et{{\text{\'et}}}
\def\an{{\text{an}}}
\def\crys{{\mathrm{crys}}}
\def\Gr{{\mathrm{Gr}}}
\def\dego{\deg=0}
\def\Qb{\bar{\Q}}
\def\ord{{\mathrm{ord}}}
\def\tr{{\mathrm{tr}}}
\def\Tr{{\mathrm{Tr}}}
\def\Frob{{\mathrm{Frob}}}
\def\cont{{\mathrm{cont}}}
\def\Ker{{\mathrm{Ker}}}    

\def\cHan#1{{\cal{H}}^{#1}_\an}
\def\cHet#1{{\cal{H}}^{#1}_\et}
\def\Han#1{{{H}^{#1}_\an}}
\def\Het#1{{{H}^{#1}_\et}}
\def\Hcont#1{{{H}^{#1}_\cont}}

\def\D{\Phi}
\def\kp{\phi_{\underline{\mu}}}
\def\cK{{\cal{K}}}
\def\cR{{\cal{R}}}

\def\qaq{\quad\text{and}\quad}
\def\qfor{\quad\text{for }}
\def\qwith{\quad\text{with }}
\def\zl{{\Z}_\ell}
\def\ql{{\Q}_\ell}
\def\qzl{\ql/\zl}
\def\zp{{\Z}_p}
\def\qp{{\Q}_p}
\def\qzp{\qp/\zp}
\def\qzl{\ql/\zl}
\def\nz{\Z/n \Z}
\def\lnz{\Z/\ell^n \Z}

\def\indlim#1{\underset{#1}{\varinjlim} \ }
\def\projlim#1{\underset{#1}{\varprojlim} \ }

\def\cXnzij{c_{\et}^{i,j}}
\def\cXnz{c_{\et}}
\def\cXlnz{c_{X,\nz}^{2,1}}
\def\cXcont{c_{\cont}}
\def\regX{\rho_{X}}

\def\rhopr{\widetilde{\rho}}
\def\rhoSpr{\widetilde{\rho}_S}
\def\rhoSpran{\widetilde{\rho}_S^\an}
\def\VSpr{\widetilde{V}_S}
\def\VSan{V^\an_S}
\def\VSpran{\widetilde{V}^\an_S}
\def\Vpr{\widetilde{V}}

\def\indlim#1{\underset{#1}{\varinjlim} \ }
\def\projlim#1{\underset{#1}{\varprojlim} \ }

\def\rmapo#1{\overset{#1}{\longrightarrow}}
\def\rmapu#1{\underset{#1}{\longrightarrow}}
\def\lmapu#1{\underset{#1}{\longleftarrow}}
\def\rmapou#1#2{\overset{#1}{\underset{#2}{\longrightarrow}}}
\def\lmapo#1{\overset{#1}{\longleftarrow}}
\def\isom{\overset{\cong}{\longrightarrow}}

\title{Chow group of 0-cycles on surface over a $p$-adic field with 
infinite torsion subgroup}
\author{Masanori Asakura and Shuji Saito}
\date\empty

\maketitle

\tableofcontents



\bigskip
\begin{abstract}
We give an example of a projective smooth surface $X$ over a $p$-adic field $K$
such that for any prime $\ell$ different from $p$, 
the $\ell$-primary torsion subgroup of $\CH_0(X)$, the Chow group of 
$0$-cycles on $X$, is infinite. 
A key step in the proof is disproving a variant of the Block-Kato conjecture 
which characterizes the image of an $\ell$-adic regulator map
from a higher Chow group to a continuous \'etale cohomology of $X$ 
by using $p$-adic Hodge theory. By aid of theory of mixed Hodge modules,
we reduce the problem to showing the exactness of de Rham complex
associated to a variation of Hodge structure, which follows
from Nori's connectivity theorem.
Another key ingredient is the injectivity result on 
\'etale cycle class map for Chow group of $1$-cycles on a proper smooth model 
of $X$ over the ring of integers in $K$ of due to K. Sato and the 
second author.
\end{abstract}
\par
\bigskip

\section{Introduction}

\medbreak

Let $X$ be a smooth projective variety over a base field $K$
and let $\CH^m(X)$ be the Chow group of algebraic cycles of 
codimension $m$ on $X$ modulo rational equivalence.
In case $K$ is a number field, there is a folklore conjecture
that $\CH^m(X)$ is finitely generated, which in particular implies
that its torsion part $\CH^m(X)_{\tor}$ is finite. 
The finiteness question has been intensively studied by many authors,
particularly for the case $m=2$ and $m=\dim(X)$ (see nice surveys
\cite{otsubo} and \cite{CT}).
\medbreak

When $K$ is a $p$-adic field (namely the completion of a number field
at a finite place), Rosenschon and Srinivas \cite{RS} have constructed the 
first example that $\CH^m(X)_{\tor}$ is infinite.
They proves that there exists a smooth projective fourfold $X$ over a 
$p$-adic field such that the $\ell$-torsion subgroup $\CH_1(X)[\ell]$ 
(see Notation) of $\CH_1(X)$, the Chow group of $1$-cycles on $X$, is 
infinite for each $\ell\in \{5,7,11,13,17\}$.
\medbreak

A purpose of this paper is to give an example of a projective smooth surface 
$X$ over a $p$-adic field such that for any prime $\ell$ different from $p$, 
the $\ell$-primary torsion subgroup $\CH_0(X)\{\ell\}$ (see Notation) of 
$\CH_0(X)$, the Chow group of $0$-cycles on $X$, is infinite. Here we note 
that for $X$ as above, $\CH_0(X)\{\ell\}$ is known to be always cofinite type 
over $\zl$  (namely the direct sum of a finite group and a finite number of 
copies of $\qzl$. The fact follows from Bloch's exact sequence 
\eqref{bloch-exact}). Thus our example presents infinite phenomena of 
different nature from the example in \cite{RS}. 
Another noteworthy point is that the phenomena discovered by our example
happens rather \it generically. \rm
\medbreak

To make it more precise, we prepare a notion of `generic surfaces' in $\P^3$.
Let 
$$
M \subset \P(H^0(\P^3_\Q,\O_\P(d)))\cong \P^{(d+3)(d+2)(d+1)/6-1}_\Q
$$
be the moduli space over $\Q$ of the nonsingular surfaces in $\P^3_\Q$
(the subscription `$\Q$' indicates the base field), and let
$$
f:{\cal X} \lra M
$$
be the universal family over $M$.
For $X \subset \P^3_K$, a nonsingular surface of degree $d$ defined over 
a field $K$ of characteristic zero, there is a morphism $t:\Spec K\to M$ 
such that $X \cong {\cal X} \times_{M} \Spec K$.
We call $X$ {\it generic} if $t$ is dominant
(i.e. $t$ factors through the generic point of $M$).
In other words, $X$ is generic if it is defined by an equation
$$
F=\underset{I}{\sum} a_{I} z^I\quad (a_I\in K)
$$ 
($[z_0:z_1:z_2:z_3]$ is the homogeneous coordinate of $\P^3$, 
$I=(i_0,\cdots,i_3)$ are multi-indices and $z^I=
z_0^{i_0}\cdots z_3^{i_3}$) satisfying the following condition:
\begin{enumerate}
\item[$(*)$]
$a_I\not=0$ for $\forall I$ and $\{a_I/a_{I_0}\}_{I\not=I_0}$
are algebraically independent over $\Q$ where $I_0=(1,0,0,0)$. 
\end{enumerate}

The main theorem is:

\begin{thm}\label{main}
Let $K$ be a finite extension of $\Q_p$
and $X\subset \P^3_K$ a nonsingular surface of degree $d\geq 5$.
Suppose that $X$ is generic and has a projective smooth model
$X_{\O_K}\subset \P^3_{\O_K}$ over the ring $\O_K$ of integers in $K$.
Let $r$ be the rank of the N\'eron-Severi group of the smooth special fiber.
Then we have
$$
\CH_0(X)\{\l\}\cong (\Q_\l/\Z_\l)^{\op r-1} \oplus \text{(finite group)}
$$
for $\l\not=p$.
\end{thm}

For example let $a_I\in K^*$ be elements such that 
$\{a_I\}_I$ are algebraically independent over $\Q$ and $\ord_K(a_I)>0$.
Then the surface $X\subset \P^3_K$ defined by an equation
$$
z_0^d-z_1^d+z_2^d-z_3^d+\underset{I}{\sum} a_I z^I\quad\text{with } d\geq 5
$$
satisfies the assumption in Theorem \ref{main} and $r\geq 2$ and hence
$\CH_0(X)$ has an infinite torsion subgroup.
Theorem \ref{main} may be compared with the finiteness results 
\cite{CTR} and \cite{R} on $\CH_0(X)_{\tor}$ for a surface $X$ over a $p$-adic 
field under the assumption that $H^2(X,\O_X)=0$ or that the rank of the 
N\'eron-Severi group does not change by reduction.

\medbreak

A distinguished role is played in the proof of Theorem \ref{main} by 
the $\ell$-adic regulator map
$$
\regX \;:\; \CH^2(X,1)\ot\Q_\l\lra 
H^1_\cont(\Spec(K),H^2(X_{\overline{K}},\Q_\l(2)))
\quad (X_{\overline{K}}=X\times_K \overline{K})
$$
from higher Chow group to continuous \'etale cohomology (\cite{cont}), 
where $\overline{K}$ is an algebraic closure of $K$ and $\ell$ is a prime
different from $\ch(K)$. It is known that the image of $\regX$ is contained 
in the subspace
$$
H_g^1(\Spec(K),V)\subset H^1_\cont(\Spec(K),V)\quad 
(V=H^2(X_{\overline{K}},\Q_\l(2)))
$$
introduced by Bloch and Kato \cite{bloch-kato}. In case $\l\not=p$ this
is obvious since $H^1_g=H^1$ by definition. For $\l=p$ this is a consequence
of a fundamental result in $p$-adic Hodge theory, which affirms that 
every representation of $G_K=\Gal(\ol{K}/K)$ arising from cohomology of 
a variety over $K$ is a de Rham representation (see a discussion
after \cite{bloch-kato}, (3.7.4)).
\medbreak

When $K$ is a number field or a $p$-adic field, 
it is proved in \cite{sato-saito2} that $\CH^2(X)\{\l\}$ is finite in case 
the image of $\regX$ coincides with $H_g^1(\Spec(K),V)$.
Bloch and Kato conjecture that it should be always the case if $K$ is 
a number field.
\medbreak

The first key step in the proof of Theorem \ref{main} is to disprove
the variant of the Block-Kato conjecture for a generic surface 
$X\subset \P^3_K$ over a $p$-adic field $K$ 
(see Theorem \ref{main2cor}). In terms of 
Galois representations of $G_K=\Gal(\ol{K}/K)$, 
our result implies the existence of a $1$-extension of 
$\Q_\l$-vector spaces with continuous $G_K$-action:
$$
0 \to H^2(X_{\overline{K}},\Q_\l(2)) \to E \to \Q_l \to 0,
\leqno(*)
$$
such that $E$ is a de Rham representation of $G_K$ 
but that there is no $1$-extension of motives over $K$:
$$
0\to h^2(X)(2) \to M \to h(\Spec(K)) \to 0
$$
which gives rise to $(*)$ under the realization functor.
The main ingredient in the proof of the first key result is 
Nori's connectivity theorem (cf. \cite{green}). 
This is done in \S3 after in \S2 we review some basic 
facts on cycle class map for higher Chow groups.
\medbreak

Another key ingredient is the injectivity result on 
\'etale cycle class map for Chow group of $1$-cycles on a proper smooth model 
of $X$ over the ring $O_K$ of integers in $K$ due to Sato and the second author
\cite{sato-saito}. It plays an essential role in deducing 
the main theorem 
\ref{main} from the first key result, which is done in \S4.
\medbreak

Finally, in \S5 Appendix, we will apply our method to produce an example of 
a curve $C$ over a $p$-adic field such that
$SK_1(C)_{\tor}$ is infinite.

\bigskip

\noindent{\it Acknowledgment}.
The authors are grateful to Dr. Kanetomo Sato for the stimulating
discussions and helpful comments.

\noindent{\it Notations}.
For an abelian group $M$,
we denote by $M[n]$ (resp. $M/n$)
the kernel (resp. cokernel) of multiplication by $n$.
For a prime number $\l$ we put
$$
M\{\l\}:=\bigcup_n M[\l^n], \quad
M_\tor:=\bigoplus_\l M\{\l\}.
$$
For a nonsingular variety $X$ over a field
$\CH^j(X,i)$ denotes Bloch's higher Chow groups.
We write $\CH^j(X):=\CH^j(X,0)$ the (usual) Chow groups. 

\bigskip

\section{Review on cycle class map and $\ell$-adic regulator}

In this section $X$ denotes a smooth variety over a field $K$ 
and $n$ denotes a positive integer prime to $\ch(K)$.

\medbreak\noindent
{\bf 2.1:}  
By \cite{GL} we have the cycle class map
$$
\cXnzij \;:\; \CH^i(X,j,\nz) \to H_\et^{2i-j}(X,\nz(i)),
$$
where the right hand side is the \'etale cohomology of $X$ with 
coefficient $\mu_{n}^{\ot i}$, Tate twist 
of the sheaf of $n$-th roots of unity. 
The left hand side is Bloch's higher Chow group with finite coefficient
which fits into the exact sequence
\begin{equation}\label{es1}
0\to \CH^i(X,j)/n \to \CH^i(X,j,\nz) \to \CH^i(X,j-1)[n] \to 0.
\end{equation}
In this paper we are only concerned with the map
\begin{equation}\label{glcyclemap}
\cXnz=c_\et^{2,1} \;:\; CH^2(X,1,\nz) \to H_\et^3(X,\nz(2)).
\end{equation}
By \cite{bloch-ogus} it is injective and its image is equal to  
$$
NH^3_\et(X,\nz(2))=\Ker\big(H^3_\et(X,\nz(2))\to
H^3_\et(\Spec(K(X)),\nz(2)),
$$
where $K(X)$ is the function field of $X$. 
In view of \eqref{es1} it implies an exact sequence
\begin{equation}\label{bloch-exact}
0\lra \CH^2(X,1)/n \rmapo{c_\et} 
NH^3_\et(X,\nz(2))\lra \CH^2(X)[n] \lra 0.
\end{equation}


\medbreak\noindent
{\bf 2.2:}  
We also need the cycle map to continuous \'etale cohomology group
(cf. \cite{cont}):
$$
\cXcont\;:\; \CH^2(X,1) \lra H^3_\cont(X,\zl(2)),
$$
where $\ell$ is a prime different from $\ch(K)$.
Note that in case $K$ is a $p$-adic field we have
$$
H^3_\cont(X,\zl(2))=\projlim n H^3_\et(X,\lnz(2))
$$
and $\cXcont$ is induced by $c_\et$ by passing to the limit.
We have the Hochschild-Serre spectral sequence
\begin{equation}\label{HSss}
E_2^{i,j}= H^i_{\cont}(\Spec(K),H^j(X_{\ol{K}},\zl(2))) \Rightarrow
H^{i+j}_\cont(X,\zl(2)).
\end{equation}
If $K$ is finitely generated over the prime subfield and $X$ is proper smooth
over $K$, the Weil conjecture proved by Deligne implies that 
$H^0(\Spec(K),H^3(X_{\ol{K}},\ql(2)))=0$.
The same conclusion holds if $K$ is a $p$-adic field and $X$ is 
proper smooth having good reduction over $K$.
Thus we get under these assumptions the following map
\begin{align}\label{highercycle1}
\rho_X \;:\; \CH^2(X ,1) \lra H^1_\cont(\Spec(K),H^2(X_{\ol{K}},\Q_\l(2)))
\end{align}
as the composite of $\cXcont\otimes\ql$ and an edge homomorphism
$$
H^{3}_\cont(X,\ql(2))\to H^1_\cont(\Spec(K),H^2(X_{\ol{K}},\Q_\l(2))).
$$

\medbreak\noindent
{\bf 2.3:}  
For later use, we need an alternative definition of cycle class maps.
For an integer $i\geq1$ we denote by $\K_i$ the
sheaf on $X_{\Zar}$, the Zariski site on $X$, associated to a presheaf
$U\longmapsto K_i(U)$.
By \cite{lands}, 2.5, we have canonical isomorphisms
\begin{equation}\label{lands}
\CH^2(X,1) \simeq H^1_{\Zar}(X,\K_2),\quad
\CH^2(X,1,\nz) \simeq H^1_{\Zar}(X,\K_2/n).
\end{equation}
Let $\epsilon^{\et}:X_\et\to X_\Zar$ be the natural map of sites and put
$$
\cH^i_\et(\nz(r))=R^i\epsilon^{\et}_* \mu^{\ot r}_n.
$$
The universal Chern classes in the cohomology groups
of the simplicial classifying space for $\GL_n$ ($n\geq 1$)
give rise to higher Chern class maps on algebraic $K$-theory
(cf. \cite{gillet}, \cite{schneider}). It gives rise to a map of sheaves:
\begin{equation}\label{et-regulator0}
\K_i/n \lra \cH_\et^i(\nz(i)).
\end{equation}
By \cite{ms} it is an isomorphism for $i=2$ and induces an isomorphism
\begin{equation}\label{et-regulator1}
H_\Zar^1(X,\K_2/n)\os{\cong}{\lra} H^{1}_\Zar(X,\cH_\et^2(\nz(2))).
\end{equation}
By the spectral sequence
$$
E_2^{pq}=H^{p}_\Zar(X,\cH_\et^q(\nz(2)))\Longrightarrow
H^{p+q}_\et(X,\nz(2)).
$$
together with the fact $H^{p}_\Zar(X,\cH_\et^q(\nz(2)))=0$ for $p>q$
shown by Bloch-Ogus \cite{bloch-ogus}, we get an injective map
\begin{equation*}\label{et-regulator3}
H^{1}_\Zar(X,\cH_\et^2(\nz(2)))\lra H^3_\et(X,\nz(2))
\end{equation*}
Again by the Bloch-Ogus theory the image of the above map
coincides with the coniveau filtration $NH^3_\et(X,\nz(2))$.
Combined with \eqref{lands} and \eqref{et-regulator1} we thus get the map
\begin{equation*}\label{et-regulator2}
c_\et:\CH^2(X,1,\nz) \os{\cong}{\lra}
H_\Zar^1(X,\K_2/n)\os{\cong}{\lra} NH^3_\et(X,\nz(2))
\os{\subset}{\lra} H^3_\et(X,\nz(2)).
\end{equation*}
One can check the map agrees with the map \eqref{glcyclemap}.

\medbreak\noindent
{\bf 2.4:}  
Now we work over the base field $K=\C$.
Let $X_\an$ be the site on the underlying analytic space $X(\C)$ endowed 
with the ordinary topology.
Let $\epsilon^{\an}:X_\an\to X_\Zar$ be the natural map of sites and put
$$
\Han i(\Z(r))=R^i\epsilon^{\an}_* \Z(r) \quad (\Z(r)=(2\pi\sqrt{-1})^r \Z).
$$
Higher Chern class map then gives a map of sheaves
\begin{equation*}\label{b-regulator0}
\K_i \lra \cHan i(\Z(i)).
\end{equation*}
By the same argument as before, it induces a map
\begin{equation}\label{b-regulator1}
c_\an:\CH^2(X,1) \os{\cong}{\lra}
H_\Zar^1(X,\K_2) \lra \Han 3(X(\C),\Z(2))
\end{equation}

\begin{lem}\label{Hodgefilt}
The image of $c_\an$ is contained in $F^2\Han 3(X(\C),\C)$, 
the Hodge filtration defined in \cite{hodge2}. 
\end{lem}

\begin{pf}
Let $\cH_D^r(\Z(i))$ be the sheaf on $X_\Zar$ associated to a presheaf 
$$
U\mapsto H^r_D(U,\Z(i))
$$ 
where $H_D^\bullet$ denotes Deligne-Beilinson cohomology 
(cf. \cite{EV}, 2.9). Higher Chern class maps to 
Deligne-Beilinson cohomology give rise to the map
$K_2 \to \cH_D^2(\Z(2))$ and $c_\an$ factors as in the following commutative 
diagram
$$
\begin{CD}
H^1_\Zar(X,\K_2)@>>>
H^1_\Zar(X,\cH_D^2(\Z(2)))@>>>
H^1_\Zar(X,\cH_\an^2(\Z(2)))\\
@.@V{a}VV@VVV\\
@.H^3_D(X,\Z(2))@>{b}>>H^3_\an(X(\C),\Z(2)).
\end{CD}
$$
Here the map $a$ is induced from the spectral sequence
$$
E_2^{pq}=H^{p}_\Zar(X,\cH_D^q(\Z(2)))\Longrightarrow
H^{p+q}_D(X,\Z(2))
$$
in view of the fact $H^{p}_\Zar(X,\cH_D^1(2))=0$ for $\forall p>0$ since
$\cH^1_D(2)\cong \C/\Z(2)$ (constant sheaf).
Since the image of $b$ is contained in $F^2\Han 3(X(\C),\C)$ 
(see \cite{EV}, 2.10), so is the image of $c_\an$.
\end{pf}
\medbreak

\begin{lem}\label{compatibility1}
We have the following commutativity diagram
\begin{equation}\label{cd1}
\begin{CD}
\CH^2(X,1) @>{c_\an}>> \Han 3(X(\C),\Z(2)) \\
@VVV @VVV \\
\CH^2(X,1,\nz) @>{c_\et}>> \Het 3(X,\nz(2)) 
\end{CD}
\end{equation}
Here the right vertical map is the composite
$$
\Han 3(X(\C),\Z(2))\to \Han 3(X(\C),\Z(2)\otimes \nz) 
\os{\cong}{\lra} \Het 3(X,\nz(2)) 
$$
and the isomorphism comes from the comparison isomorphism
between \'etale cohomology and ordinary cohomology 
(SGA$4\frac{1}{2}$, Arcata, 3.5) together with the isomorphism
$$
\Z(1)\otimes\nz \simeq (\epsilon^\an)^* \mu_n
$$ 
given by the exponential map.
\end{lem}

\begin{pf}
This follows from the compatibility of the universal Chern classes
(\cite{gillet} and \cite{schneider}).
\end{pf}

\bigskip

\section{Counterexample to Bloch-Kato conjecture over $p$-adic field}

In this section $K$ denotes a $p$-adic field and let $X$ be a proper smooth
surface over $K$. We fix a prime $\ell$ (possibly $\ell=p$) and 
consider the map \eqref{highercycle1}
\begin{equation}\label{highercycle2}
\rho_X \;:\; \CH^2(X,1) \lra 
H^1_\cont(\Spec(K),V) \quad (V=H^2_\et(X_{\ol{K}},\Q_\l(2)))).
\end{equation}
Define the primitive part $\Vpr$ of $V$ by:
\begin{equation}\label{Vpr}
\Vpr:= H^2_\et(X_{\ol{K}},\Q_\l(2))/V_0, \quad V_0=[H_X]\ot\Q_\l(1),
\end{equation}
where $[H_X]\in H^2_\cont(X_{\ol{K}},\Q_\l(1))$
is the cohomology class of a hyperplane section. Let
$$
\rhopr \;:\; \CH^2(X,1) \lra H^1_\cont(\Spec(K),\Vpr)
$$
be the induced map.

\begin{thm}\label{main2}
Let $X\subset \P^3_K$ be a generic smooth surface of degree $d\geq 5$. 
Then $\rhopr$ is the zero map for arbitrary $\l$.
\end{thm}

\begin{rem}\label{rem1}
\begin{itemize}
\item[(1)] 
The key point of the proof is Nori's connectivity theorem.
This is an analogue of \cite{voisin} 1.6 (where she worked
on Deligne-Beilinson cohomology).
\item[(2)]
Bloch-Kato \cite{bloch-kato} considers regulator maps such as 
\eqref{highercycle2} for a smooth projective variety over a number 
field and conjectures that its image coincides with $H^1_g$.
We will see later (see Theorem \ref{main2cor}) 
that the variant of the conjecture over 
a $p$-adic field is false in general. 
\item[(3)]
The construction of a counterexample mentioned in (2) hinges on the
assumption that the surface $X\subset\P^3_K$ is generic. One may still
ask whether the image of $l$-adic regulator map coincides $H^1_g$
for a proper smooth variety $X$ over a $p$-adic field when $X$ is defined
over a number field.
\end{itemize}
\end{rem}
\medbreak

\begin{pf}
Let $f:\cX \to M$ be as in the introduction and let $t:\Spec(K) \to M$
be a dominant morphism such that $X\simeq \cX\times_M \Spec(K)$.
For a morphism $S \to M$ of smooth schemes over $\Q$ let
$f_S: X_S=\cX\times_M S \to S$ be the base change of $f$.
The same construction of \eqref{highercycle1} give rise to the regulator map
$$
\rho_S\;:\; \CH^2(X_S,1) \to H^1_\cont(S,V_S),
$$
where $V_S=R^2(f_S)_*\Q_l(2)$ is a smooth $\Q_l$-sheaf on $S$. 
Define the primitive part of $V_S$:
$$
\VSpr= R^2(f_S)_*\Q_l(2)/[H]\otimes\Q_l(1),
$$
where
$[H]\in H^0(S,R^2(f_S)_*\Q_l(1))$ is the class of a hyperplane section.
Let
$$
\rhoSpr \;:\; \CH^2(X_S,1) \to H^1_\cont(S,\VSpr),
$$
be the induced map. Note
$$
\CH^2(X,1) =\indlim {S} \CH^2(X_S,1),
$$
where $S \to M$ ranges over the smooth morphisms which factor 
$t:\Spec(K)\to M$. Note also that we have the commutative diagram for such $S$:
$$
\begin{CD}
\CH^2(X_S,1) @>{\rhoSpr}>> H^1_\cont(S,\VSpr)\\
@VVV @VVV\\
\CH^2(X,1) @>{\rhopr}>> H^1_\cont(\Spec(K),\Vpr).\\
\end{CD}
$$
Thus it suffices to show
$$
H^1_\cont(S,\VSpr)=0.
$$
Without loss of generality we suppose $S$ is a affine smooth variety
over a finite extension $L$ of $\Q$.

\begin{claim}\label{claim1}
Assume $d\geq 4$. The natural map
$$
H^1_\cont(S,\VSpr) \lra H^1_\et(S_{\ol{\Q}},\VSpr)
\quad (S_{\ol{\Q}}:=S\times_L \Spec(\ol{\Q}))
$$
is injective. 
\end{claim}

Indeed, by the Hochschild-Serre spectral sequence, it is enough to see
$H^0_\et(S_{\ol{\Q}},\VSpr)=0$, which follows from \cite{AS}, Th.5.3(2).
\medbreak

By SGA$4\frac{1}{2}$, Arcata, Cor.(3.3) and (3.5.1) we have 
$$
H^1_\et(S_{\ol{\Q}},\VSpr)\cong H^1_\et(S_\C,\VSpr) \simeq 
H^1_\an(S(\C),\VSpran)\otimes\Q_l,
\quad (S_{\C}:=S\times_L \Spec(\C))
$$
where
$\VSpran$ is the primitive part of $\VSan=R^2(f_S^\an)_* \Q(2)$ 
with $f_S^\an:(X_{S_\C})_\an \to (S_\C)_\an$, the natural map of sites.
By definition $\VSpran$ is a local system on $S(\C)$ whose fiber over 
$s\in S(\C)$ is the primitive part of $H^2_\an(X_s(\C),\Q(2))$ 
for $X_s$, the fiber of $X_S\to S$ over $s$.
Due to Lemma \ref{compatibility1}, it suffices to show that the triviality
of the image of the map 
$$
\rhoSpran \;:\; \CH^2(X_{S_\C},1) \lra H^1_\an(S(\C),\VSpran)
$$
which is induced from 
$$
c_\an:  \CH^2(X_{S_\C},1) \lra H^3_\an(X_S(\C),\Q(2))
$$
by using the natural map
$$
H^3_\an(X_S(\C),\Q(2)) \to H^1_\an(S(\C),\VSan)
$$
arising from the Leray spectral sequence for 
$f_S^\an:(X_{S_\C})_\an \to (S_\C)_\an$ and the vanishing 
$R^3(f_S^\an)_*\Q(2)=0$.

\begin{claim}\label{claim2}
The image of $\rhoSpran$ is contained in the Hodge filtration
$F^2H^1_\an(S(\C),\VSpran\otimes\C)$ defined by theory of Hodge modules 
\cite{msaito}.
\end{claim}

This follows from the functoriality of Hodge filtrations and 
Lemma \ref{Hodgefilt}.
\medbreak

\begin{claim}\label{claim3}
For integers $m,p\geq 0$ there is a natural injective map
$$F^p H^m_\an(S(\C),\VSpran\otimes\C) \to H^m(S_\C,G^p DR(\VSpran))$$
where 
$G^p DR(\VSpran)$ is the complex of sheaves on $S_\C$:
\begin{multline*}
F^p H^2_\dR(X_S/S)_\prim \otimes\O_{S_\C} \to
F^{p-1} H^2_\dR(X_S/S)_\prim\ot\Omega^1_{S_\C/\C}\to \cdots\\
\cdots \to F^{p-r} H^2_\dR(X_S/S)_\prim\ot\Omega^r_{S_\C/\C}
\to F^{p-r} H^2_\dR(X_S/S)_\prim\ot\Omega^{r+1}_{S_\C/\C}\to \cdots \\
\end{multline*}
Here $H^\bullet_\dR(X_S/S)$ denotes the de Rham cohomology of $X_S/S$, and 
$H^\bullet_\dR(X_S/S)_\prim$ is its primitive part defined by the same way as 
before, and the maps are induced from the Gauss-Manin connection.
\end{claim}

This follows from \cite{asakura} Lemma 4.2.
We note that its proof hinges on theory of mixed Hodge modules. A key point
is degeneration of Hodge spectral sequence for cohomology with coefficient.
\medbreak

By the above claims we are reduced to show the exactness at the middle term
of the following complex:
\begin{equation}\label{lemma1-2}
F^2H^2_\dR(X_S/S)_\prim \otimes\O_{S_\C} \to
F^1H^2_\dR(X_S/S)_\prim\ot\Omega^1_{S_\C/\C}\to
H^2_\dR(X_S/S)_\prim\ot\Omega^2_{S_\C/\C}.
\end{equation}
For this it suffices to show that
$$
f_*\Omega^2_{X_S/S}\otimes\O_{S_\C}\to
(R^1f_*\Omega^1_{X_S/S})_\prim\ot\Omega^1_{S_\C/\C}\to
R^2f_*\O_{X_S}\ot\Omega^2_{S_\C/\C}
$$
is exact at the middle term and 
$$
f_*\Omega^2_{X_S/S}\ot\Omega^1_{S_\C/\C}\to
(R^1f_*\Omega^1_{X_S/S})_\prim\ot\Omega^2_{S_\C/\C}
$$
is injective.
The former holds when $d\geq 5$
and the latter holds
when $d\geq 3$
by Nori's connectivity theorem (cf. \cite{green}).
This completes the proof of \ref{main2}.
\end{pf}

\bigskip

Let $\O_K\subset K$ be the ring of integers and $k$ be the residue field. 
In order to construct an example where the image of the regulator map 
$$
\rho_X \;:\; \CH^2(X,1) \rmapo{\rho_X} H^1_\cont(\Spec(K),V) \quad
(V=H^2_\et(X_{\ol{K}},\Q_\l(2))))
$$
is not equal to $H_g^1(\Spec(K),V)$, we now take a proper smooth surface 
$X$ having good reduction over $K$ so that $X$ has a proper smooth model 
$X_{\O_K}$ over $\Spec(\O_K)$. We denote the special fiber by $Y$.
By \cite{langer-saito} (see the diagram below 5.7 on p.341), there is a 
commutative diagram
\begin{equation}\label{LSCD}
\begin{CD}
\CH^2(X,1) @>{\rhopr}>> H^1_g(\Spec(K),V) \\
@VV{\partial}V @VVV \\
\CH^1(Y) @>{\alpha}>> H^1_\cont(\Spec(K),V)/H^1_f(\Spec(K),V)\\
\end{CD}
\end{equation}
where $H^1_f\subset H^1_g\subset H^1_\cont$
are the subspaces introduced by Bloch-Kato \cite{bloch-kato} and 
$\partial$ is a boundary map in localization sequence for higher Chow groups.

\begin{thm}\label{main2cor}
Let $X\subset \P^3_K$ be a generic smooth surface of degree $d\geq 5$. 
Assume that $X$ has a projective smooth model $X_{\O_K}\subset \P^3_{\O_K}$
over $\O_K$ and let $Y\subset \P^3_k$ be its special fiber.
\begin{itemize}
\item[(1)]
The image of $\partial\otimes\Q$ is contained in the subspace of 
$\CH^1(Y)\otimes\Q$ 
generated by the class $[H_Y]$ of a hyperplane section of 
$Y$.
\item[(2)]
Let $r$ be the Picard number of $Y$. Then 
$$
\dim_{\Q_\l}\big(H^1_g(\Spec(K),V)/\Image(\rho_X)\big)\geq r-1.
$$
\end{itemize}
\end{thm}
\begin{pf}
Letting $V_0\subset V$ be as in \eqref{Vpr}, we have a decomposition 
$V=\Vpr\oplus V_0$ as $G_K$-modules. Let $W\subset CH^2(X,1)$ be the image of 
$\Z\cdot [H_X] \otimes K^\times$ under the product map
$\CH^1(X)\ot K^\times \to \CH^2(X,1).$
Then it is easy to see $\regX$ induces an isomorphism
$$
W\ot\ql \simeq H^1_g(\Spec(K),V_0)=H^1_\cont(\Spec(K),V_0)
$$
and that $\partial(W)=\Z\cdot[H_Y]\subset \CH^1(Y)$.
Hence (1) 
follows from Theorem \ref{main2} together with injectivity of 
$\alpha$ in \eqref{LSCD} proved by \cite{langer-saito}, Lemma 5-7.

As for (2) we first note that
$\dim_{\ql}\big(H^1_\cont(\Spec(K),V_0)/H^1_f(\Spec(K),V_0)\big)=1$ 
(see \cite{bloch-kato}, 3.9).
Moreover the same argument (except using the Tate conjecture) in the last part 
of \S5 of \cite{langer-saito} shows 
$$
\dim_{\ql}(\CH^1(Y)\ot\ql) \leq 
\dim_{\ql}\big(H^1_g(\Spec(K),V)/H^1_f(\Spec(K),V)\big).
$$
Hence (2) follows from (1).
\end{pf}

\medbreak

\begin{rem}
Let the assumption be as in Corollary \ref{main2cor}. Then we have 
$$
\dim_{\Q_\l}\left(H^1_g(\Spec(K),V)/\Image(\rho_X)\right)\geq
\begin{cases}
r-1& \l\not=p\\
r-1+(h^{0,2}+h^{1,1}-1)[K:\Q_p]& \l=p
\end{cases}
$$
where $h^{p,q}:=\dim_K H^q(X,\Omega^p_{X/K})$ denotes the Hodge number.
Moreover the equality holds if and only if the Tate conjecture for
divisors on $Y$ holds. It follows from Theorem \ref{main2} and computation of
$\dim_{\Q_\l}H^1_g(\Spec(K),V)$ using \cite{bloch-kato} 3.8 and 3.8.4.
The details are omitted.
\end{rem}

\bigskip

\section{Proof of Theorem \ref{main}}

Let $K$ be a $p$-adic field and $\O_K\subset K$ the ring of 
integers and $k$ the residue field. Let us consider schemes
\begin{equation}\label{setup1}
\begin{CD}
X@>{j}>> X_{\O_K} @<<{i}< Y\\
@VVV@VVV@VVV\\
\Spec(K) @>>> \Spec(\O_K) @<<< \Spec(k)
\end{CD}
\end{equation}
where all vertical arrows are projective and smooth of relative dimension $2$
and the diagrams are Cartesian.
We have a boundary map in localization sequence for higher Chow groups
with finite coefficient
$$
\partial\;:\; \CH^2(X,1,\nz) \to \CH^1(Y)/n.
$$
For a prime number $\l$, it induces 
$$
\partial_\l \;:\; \CH^2(X,1,\qzl) \to \CH^1(Y)\otimes\qzl,
$$
where
$\CH^2(X,1,\qzl):=\indlim n \CH^2(X,1,\lnz)$.

\begin{thm}\label{main3}
For $\l\not=p:=\ch(k)$, $\partial_\l$ is surjective and finite kernel.
Hence we have
$$
\CH^2(X,1,\qzl) 
\cong (\qzl)^{\op r}+\text{(finite group)}
$$
where $r$ is the rank of $\CH^1(Y)$.
\end{thm}

Theorem \ref{main} is an immediate consequence of 
Theorems \ref{main2cor} (1), 
\ref{main3}, and the exact sequence \eqref{es1}:
$$
0\to \CH^2(X,1)\otimes\qzl  \to \CH^2(X,1,\qzl) \to \CH^2(X)\{\l\} \to 0.
$$

\begin{pf}
Write $\Lambda=\qzl$
We have a commutative diagram:
\begin{equation*}
\begin{CD}
\CH^2(X,1,\Lambda) @>{\partial}>> 
\CH^1(Y)\otimes\Lambda  @>{i_*}>>
\CH^2(X_{\O_K})\otimes\Lambda @>{j^*}>> 
\CH^2(X)\otimes\Lambda \\
@VV{c_1}V @VV{c_2}V @VV{c_3}V @VV{c_4}V \\
\Het 3(X,\Lambda(2)) @>{\partial_\et}>> \Het 2(Y,\Lambda(1))  @>{i_*}>>
\Het 4(X_{\O_K},\Lambda(2)) @>{j^*_\et}>> \Het 4(X,\Lambda(2)) 
\end{CD}
\end{equation*}
Here the upper (resp. lower) exact sequence arises from localization theory 
for higher Chow groups with finite coefficient (resp. \'etale cohomology
together with absolute purity \cite{purity}). The vertical maps are \'etale 
cycle class maps. By \eqref{bloch-exact} $c_1$ is injective. Noting 
$\CH^1(Y)\simeq \Het 2(Y,\G_m)$, $c_2$ is injective by the Kummer theory. 
It is shown in \cite{sato-saito} that $c_3$ is an isomorphism.
Hence the diagram reduces the proof of Theorem \ref{main3} 
to showing that 
$\Ker(\partial_\et)$ and $\Ker(j^*_\et)$ are finite.
This is an easy consequence of the proper base change theorem for 
\'etale cohomology and the Weil conjecture (\cite{weilconj}). 
For the former we use also an exact sequence
$$
\Het 3(X_{\O_K},\Lambda(2)) \to \Het 3(X,\Lambda(2)) \rmapo{\partial_\et} 
\Het 2(Y,\Lambda(1)).
$$
\end{pf}

\bigskip

\section{Appendix: $SK_1$ of curve over $p$-adic field}

Let $C$ be a proper smooth curves over a field $K$ and consider
$\CH^2(C,1)$. By \cite{lands}, 2.5, we have an isomorphism
$$
\CH^2(C,1) \simeq H^1_\Zar(C,\K_2) \simeq SK_1(C).
$$
By definition
$$
SK_1(C) =\Coker(K_2(K(C)) \rmapo{\delta} 
\underset{x\in C_0}{\bigoplus} K(x)^\times),
$$
where
$K(C)$ is the function field of $C$, $C_0$ is the set of the closed points of
$C$, and $K(x)$ is the residue field of $x\in C_0$, and $\delta$ is given
by the tame symbols. The norm maps
$K(x)^\times \to K^\times$ for $x\in C_0$ induce
$$
N_{C/K}\;:\; SK_1(C) \to K^\times.
$$
We write $V(C)=\Ker(N_{C/K})$.
\medbreak

When $K$ is a $p$-adic field, it is known 
by class field theory for curves over local field (\cite{saito}),
that $V(C)$ is a direct sum of its maximal divisible subgroup
and a finite group. An interesting question is whether the divisible subgroup
is uniquely divisible, or equivalently whether $SK_1(C)_{\tor}$ is finite.
In case the genus $g(C)=1$ affirmative results have been obtained in 
\cite{tsato} and \cite{asakura2}. 
The purpose of this section is to show that the method in the previous 
sections gives rise to an example of a curve $C$ of $g(C)\geq 2$ such that
$SK_1(C)_{\tor}$ is infinite.
\medbreak

Let $C$ be as in the beginning of this section and let $n$ be
a positive integer prime to $\ch(K)$. We have the cycle class map
\begin{equation}\label{clsk}
c_\et\;:\; \CH^2(C,2,\nz) \to \Het 2 (C,\nz(2)).
\end{equation}
The main result of \cite{ms} implies that the above map is an isomorphism.
In view of the exact sequence (cf. \eqref{es1}):
$$
0\to \CH^2(C,2)/n \to \CH^2(C,2,\nz) \to SK_1(C)[n] \to 0,
$$
we get the exact sequence (\cite{suslin} 23.4):
\begin{equation}\label{suslin-es}
0\to \CH^2(C,2)/n \to \Het 2 (C,\nz(2)) \to SK_1(C)[n] \to 0.
\end{equation}
We will use also cycle class map to continuous \'etale cohomology:
$$
c_\cont \;:\; \CH^2(C,2)\otimes\Q_\l \to \Hcont 2 (C,\ql(2))
$$
where $\ell$ is any prime number different from $\ch(K)$.
When $K$ is a $p$-adic field, one easily shows 
\begin{equation}\label{HSeq}
\Hcont 2 (C,\ql(2)) \simeq \Hcont 1(\Spec(K),\Het 1 (C_{\ol{K}},\ql(2)))
\end{equation}
by using the Hochschild-Serre spectral sequence \eqref{HSss}.
Hence we get the map
\begin{equation}\label{regC}
\rho_C \;:\; \CH^2(C,2)\otimes\Q_\l \to 
\Hcont 1(\Spec(K),\Het 1 (C_{\ol{K}},\ql(2))).
\end{equation}
\bigskip

Let $M_g$ be the moduli space of tri-canonically embedded projective 
nonsingular curves of genus $g\geq 2$ over the base field $\Q$ 
(cf. \cite{DM}), and let $f:{\cal C}\to M_g$ be the universal family.
For $C$ over $K$ be as before, we say $C_K$ is generic if there is a 
dominant morphism $\Spec(K) \to M_g$ such that 
$C\cong {\cal C}\times_{M_g}\Spec(K)$.

\begin{thm}\label{main4}
Let $K$ be a $p$-adic field and let $C$ be a generic curve of genus 
$g\geq 2$ over $K$. Then $\rho_C$ is the zero map. We have an isomorphism
$$
SK_1(C)_{\tor}
\cong H^2_\et(C,\Q/\Z(2)) \big(:=\indlim n H^2_\et(C,\nz(2))\big).
$$
\end{thm}
\begin{rem}
This is an analogue of \cite{green-griffiths} where they
worked on Deligne-Beilinson cohomology.
\end{rem}
\begin{pf}
The second assertion follows easily from the first in view of 
\eqref{suslin-es}.
The fist assertion is shown by the same method as the proof of 
Theorem \ref{main2}, noting the following fact: 
Nori's connectivity theorem holds for generic curve of genus
$\geq 2$, namely for any dominant smooth morphism $S\to M_g$,
letting $f:C_S:={\cal C}\times_{M_g}S\to S$,
the map
$$
f_*\Omega^1_{C_S/S}\lra R^1f_*\O_{C_S}\ot\Omega^1_{S/\Q}
$$
induced from the Gauss-Manin connection is injective.
\end{pf}

\begin{cor}\label{main4cor}
Let $C$ be as in Theorem \ref{main4}. Assume the Jacobian variety
$J(C)$ has semistable reduction over $K$.
Let $\cJ$ be the Neron model of $J$ with $\cJ_s$, its special fiber.
Let $r$ be the dimension of the maximal split torus in $\cJ_s$.
For a prime $\ell$, we have
$$
SK_1(C)\{\ell\} \simeq (\qzl)^{r_\l} \oplus  \text{(finite group)},
$$
where $r_\l =r$ for $\l\not=p$ and $r_p=r+2g[K:\Q_p]$.
\end{cor}

For example $SK_1(C)\{\ell\}$ is infinite for any $\ell$ if $C$ is a
generic Mumford curve.
\medbreak

Theorem \ref{main4cor} follows from Theorem \ref{main4} 
and the following:

\begin{lem}\label{main4lem}
Let $C$ be proper smooth curve over a $p$-adic field $K$
Assume $J(C)$ has semistable reduction over $K$ and let $r_\l$ be as above.
Then 
$$
\dim_{\ql} H^2_\cont(C,\ql(2))=\dim_{\ql} H^1_\cont(\Spec(K),V) =r_\l.
\quad \big(V=H^1_\et(C_{\ol{K}},\ql(2))\big).
$$
\end{lem}
\begin{pf}
The first equality follows from \eqref{HSeq}.
By \cite{jannsen}, Th. 5 and Cor. 7 (p. 354--355) we have
$H^0_\cont(\Spec(K),V)=0$ and $\dim_{\ql} H^2_\cont(\Spec(K),V)=r$.
Lemma \ref{main4lem} now follows from computation of Euler-Poincar\'e
characteristic given in \cite{serre}, II 5.7.
\end{pf} 

\begin{rem}
It is shown that $ H^1_\cont(\Spec(K),V)=H^1_g(\Spec(K),V)$.
Hence, if $C$ is a generic curve of genus$\geq 2$, then the map 
$\rho_C$ \eqref{regC} does not surjects on to $H^1_g$ if $r_\l\geq 1$.
This gives another counterexample to a variant of the Bloch-Kato
conjecture for $p$-adic fields.
\end{rem}

\bigskip

\noindent
Graduate School of Mathematics, Kyushu University,
Fukuoka 812-8581,
JAPAN

\medskip
\noindent
E-mail address : asakura@math.kyushu-u.ac.jp

\bigskip\noindent
Graduate School of Mathematical Sciences, Tokyo University,
Tokyo 153-8914,
JAPAN

\medskip
\noindent
E-mail address : sshuji@msb.biglobe.ne.jp

\end{document}